\documentclass{article}
\usepackage{graphicx}
\usepackage{amsmath}
\usepackage{amsfonts}
\usepackage{amssymb}

\newtheorem{theorem}{Theorem}
\newtheorem{acknowledgement}{Acknowledgement}

\newtheorem{definition}{Definition}
\newtheorem{example}{Example}

\newtheorem{lemma}{Lemma}

\newtheorem{proposition}{Proposition}
\newtheorem{remark}{Remark}

\newenvironment{proof}[1][Proof]{\textbf{#1.} }{\ \rule{0.5em}{0.5em}}

\begin{document}

\title{Generalized Casimir invariants of six dimensional solvable real Lie algebras with five dimensional nilradical }

\author{Rutwig Campoamor-Stursberg\\Laboratoire de Math\'ematiques et Applications\\F.S.T Universit\'e de Haute Alsace\\
F-68093 Mulhouse (France)\\R.Campoamor@uha.fr}

\date{}

\maketitle

\begin{abstract}
We finish the determination of the invariants of the coadjoint representation of six dimensional real Lie algebras, by determining a fundamental set of invariants for the 99 isomorphism classes of solvable Lie algebras with five dimensional nilradical. We also give some results on the invariants of solvable Lie algebras in arbitrary dimension whose nilradical has codimension one.
\end{abstract}

\section{Introduction}

Invariant functions of the coadjoint orbit of semisimple Lie algebras were determined and physically interpreted by Racah in 1950, although some precedent existed in the literature, starting with the pioneering work of Casimir in 1931 (which justifies that these functions are known as Casimir invariants) \cite{Cas}. More specifically, the polynomials in the universal enveloping algebra invariant by the adjoint representation of a Lie algebra were studied in the case of classical groups by Pontryagin, and by Yen for the exceptional groups \cite{Yen}. For semisimple Lie algebras the number of functionally independent invariants is completely determined by the dimension of the Cartan subalgebra. This fact led to the belief that other groups would provide invariants of similar form. But for solvable Lie algebras the results found were of a very different nature. Neither the number of invariants nor their specific structure can be determined generically. This fact is due to the more complicated structure of solvable Lie algebras, which do not underly to a root theory like the classical case. During the last decades, the invariants of solvable Lie algebras have become more important due to their significance in representation theory of Lie algebras, where they allow to label irreducible representations or to split representations \cite{Mi}. In physics they have also shown their importance, as for example in the construction of completely integrable hamiltonian dynamical systems or the obtention of quantum numbers for symmetry groups of systems. Invariants have been determined only for a small number of solvable and non-semisimple Lie algebras, like the Poincar\'e algebra and its subalgebras or solvable Lie algebras in dimensions less or equal to six \cite{Ge,PWZ}. Indeed the solvable case is a quite important one, mainly because of the Levi decomposition theorem, which shows that any Lie algebra is the semidirect product of a semisimple subalgebra and its radical, which is solvable. Since the invariants of semisimple algebras are perfectly known, the question reduces to the invariants of the radicals and the action of the Levi subalgebra on the radicals.\newline
Solvable Lie algebras have been classified over the real and complex numbers up to dimension six, and their invariants for algebras having rank $r\leq 4$ were computed by Ndogmo few years ago \cite{Nd}. However, this work is incomplete, since the author left out the case of six dimensional Lie algebras whose nilradical is of dimension five. The justification was that their structure is very similar to the five dimensional case with four dimensional nilradical \cite{PWZ}. In certain aspects this omission is justified, since the structure of the invariants is related, but it must be kept in mind that odd dimensional Lie algebras must have at least an invariant, while for even dimensional Lie algebras it can happen that only trivial invariants exist, i.e., that any invariant function is necessarily a constant. Exactly this point makes this case interesting and worthy to be studied in detail, and allows to obtain some results which are interesting for the general case, as well as the fact that this classification has not been revised due to the relative inaccesibility of the original reference, up to the discovery of some errors by Turkowski when finishing the classification \cite{Tu1}.\newline In this paper we determine the invariants of all real Lie algebras in dimension six whose nilradical is of dimension five, and we show how from the tables we can infer some results which can be formulated in arbitrary dimension. In particular, we show how to construct parametrized families of solvable Lie algebras whose nilradical has codimension one and which have no non-trivial invariants.

\bigskip

Unless otherwise stated, any Lie algebra $\frak{g}$ considered here is of finite dimension over the field $\mathbb{K}=\mathbb{R}$.\newline The notation used for the six dimensional Lie algebras with five dimensional nilradical, as well as the nomenclature for the latter, are extracted from the paper by Mubarakzyanov where they were classified \cite{Mu3}. We also note that these lists contain some minor mistakes, which have been corrected here without explicit mention.

\section{The coadjoint representation}

Let $G$ be a Lie group and $\frak{g}$ the corresponding Lie algebra. The
coadjoint representation of $G$ is given by the mapping:
\begin{equation*}
ad^{\ast }:G\rightarrow GL\left( \frak{g}^{\ast }\right) :\;\left(
ad_{g}^{\ast }F\right) \left( x\right) =F\left( ad_{g^{-1}}x\right) ,\;g\in
G,F\in \frak{g}^{\ast },x\in \frak{g}
\end{equation*}

\begin{definition}
A function $F\in C^{\infty }\left( \frak{g}^{\ast }\right) $ is called a
semi-invariant for the coadjoint representation if
\begin{equation*}
F\left( ad_{g^{-1}}x\right) =\chi \left( g\right) F\left( x\right)
\end{equation*}
for any $g\in G$, $\chi $ being a character of $\frak{g}$. In particular, if
$\chi \left( g\right) =1$, then $F$ is called an invariant of $\frak{g}$.
\end{definition}

Thus $F$ is an invariant if and only if $F$ is constant on each coadjoint orbit.
Semi-invariants of the coadjoint representation are of extraordinary
importance for the analysis of completely integrable Hamiltonian dynamical
systems, as pointed out by Arhangelsky in \cite{Ar}, and were used by Trofimov
to determine the invariants of Borel subalgebras of simple Lie algebras \cite
{Tro1} by means of the symmetries of root systems.

\begin{definition}
A maximal set of functionally independent invariants of $\frak{g}$ is called
a fundamental set of invariants.
\end{definition}

The usual method to compute the (semi-)invariants of a Lie algebra is making
use of the theory of linear partial differential equations \cite{PP,Mir,Fo}. Let $%
\left\{ X_{1},..,X_{n}\right\}$ be a basis of $\frak{g}$ and let $\left\{
C_{ij}^{k}\right\} $ be its structure constants over this basis. We can
represent $\frak{g}$ in the space $C^{\infty }\left( \frak{g}^{\ast }\right)
$ by the differential operators
\begin{equation}
\widehat{X}_{i}=-C_{ij}^{k}x_{k}\frac{\partial }{\partial x_{j}}
\end{equation}
where $\left[ X_{i},X_{j}\right] =C_{ij}^{k}X_{k}$ \ $\left( 1\leq i<j\leq
n\right) $. It is not difficult to verify that the operators $\widehat{X}_{i}
$ satisfy the brackets $\left[ \widehat{X}_{i},\widehat{X}_{j}\right]
=C_{ij}^{k}\widehat{X}_{k}$. The following result gives the procedure to
obtain the invariants of the coadjoint representation:

\begin{theorem}
A function $F\in C^{\infty }\left( \frak{g}^{\ast }\right) $ is an invariant
if and only if it is a solution of the following system:
\begin{equation}
\left\{\widehat{X}_{i}F=0,\;1\leq i\leq n\right\}
\end{equation}
\end{theorem}

Proofs can be found in \cite{Fo,LL,Ro,Tro2}. This theorem reduces the
determination of the invariant to a system of linear first-order partial
differential equations. The method has been used by many authors to study
the invariants, and it is the one we use here. In particular, this allows to
determine the classical Casimir operators of any Lie algebra \cite{Cas,Ra}.
If $\frak{A}$ denotes the universal enveloping algebra of $\frak{g}$ and $%
Z\left( \frak{A}\right) $ its center, then it is well known that the
elements in $Z\left( \frak{A}\right) $ correspond to the Casimir operators
\cite{Schw}. This set indeed coincides with the set of polynomial invariants
of $ad^{\ast }$, while the rational invariants correspond to ratios of
polynomials contained in $\frak{A}$. However, the system  $\left( 2\right) $
can also have non-rational solutions, which do not have such an
interpretation in terms of the enveloping algebra. Such solutions will be
called generalized Casimir invariants.

\medskip

For any given Lie algebra $\frak{g}$, the number of functionally independent invariants of the coadjoint representation can be computed from the brackets. More specifically, let $\left(C_{ij}^{k}x_{k}\right)$ be the matrix which represents the commutator table over the basis $\left\{X_{1},..,X_{n}\right\}$, the $C_{ij}^{k}$ being the structure constants over this basis. The matrix is clearly skew-symmetric, which implies that the rank is necessarily even.

\begin{proposition}
\cite{Be} The cardinal $\mathcal{N}$ of a fundamental set of
invariants of $\frak{g}$ is given by $\dim \,\frak{g}-\sup \left\{
rank\left( C_{ij}^{k}x_{k}\right) _{1\leq i<j\leq \dim \frak{g}}\right\} $.
In particular, $dim \frak{g}\equiv\mathcal{N} \left(mod\quad 2\right) $.
\end{proposition}

Thus the number of functionally independent invariants is deduced from the
dimension of $\frak{g}$ and the maximal rank of its commutation table
considered as matrix, and has the same parity as the dimension of the
algebra. In general, the number of polynomial solutions of $\left( 2\right) $
will be strictly lower than $\dim \,\frak{g}-\sup \left\{ rank\left(
C_{ij}^{k}x_{k}\right) _{1\leq i<j\leq \dim \frak{g}}\right\} $, and only
for particular classes of Lie algebras, such as semi-simple or nilpotent Lie
algebras we will obtain an equality \cite{Ra,AA,So}. On the other hand, for
algebraic Lie algebras (i.e., the Lie algebras of algebraic groups \cite{Hu}%
) all invariants can be chosen as rational functions \cite{Di}. In recent
years, and due to the impossibility of classifying solvable Lie algebras in
dimensions $n\geq 7$, the analysis has been restricted to important classes of algebras, such as solvable Lie algebras having abelian or Heisenberg
nilradicals \cite{Ru,Wi} or rigid Lie algebras \cite{Ca,Ca2}.

\medskip

We illustrate the determination of the invariants by an example:

\begin{example}
Let $\frak{g}_{6,65}$ be the six dimensional (real)  Lie algebra whose
brackets are:
\begin{equation*}
\begin{array}{l}
\left[  X_{3},X_{5}\right]  =X_{1}, \quad \left[
X_{4},X_{5}\right]  =X_{2}\\
\left[  X_{1},X_{6}\right]  =\lambda X_{1}+X_{2},\quad \left[  X_{2}%
,X_{6}\right]  =\lambda X_{2}\\
\left[  X_{3},X_{6}\right]  =\left(  \lambda-\gamma\right)  X_{3}+X_{4},\quad \left[  X_{4},X_{6}\right]  =\left(  \lambda-\gamma\right)  X_{4}\\
\left[  X_{5},X_{6}\right]  =\gamma X_{5}
\end{array}
\end{equation*}
In this case, system $\left( 2\right) $ has the form
\begin{equation}
\left.
\begin{array}{l}
\widehat{X}_{1}F=-\left( \lambda x_{1}+x_{2}\right) \partial _{x_{6}}F=0 \
\widehat{X}_{2}F=-\lambda x_{2}\partial _{x_{6}}F=0 \\
\widehat{X}_{3}F=\left( -x_{1}\partial _{x_{5}}-\left( \left( \lambda -\mu
\right) x_{3}+x_{4}\right) \partial _{x_{6}}\right) F=0 \\
\widehat{X}_{4}F=\left( -x_{2}\partial _{x_{5}}-\left( \lambda -\mu \right)
x_{4}\partial _{x_{6}}\right) F=0 \\
\widehat{X}_{5}F=\left( x_{1}\partial _{x_{3}}+x_{2}\partial _{x_{4}}-\gamma
x_{5}\partial _{x_{6}}\right) F=0 \\
\widehat{X}_{6}F=\left( \lambda x_{1}+x_{2}\right) \partial _{x_{1}}+\lambda
x_{2}\partial _{x_{2}}+\left( \left( \lambda -\mu \right) x_{3}+x_{4}\right)
\partial _{x_{3}}+ \\
+\left( \lambda -\mu \right) x_{4}\partial _{x_{4}}+\gamma x_{5}\partial
_{x_{5}}=0
\end{array}
\right\}
\end{equation}
The system reduces to the following
\begin{equation}
:\left.
\begin{array}{l}
\widehat{X}_{5}^{\prime }F=\left( x_{1}\partial _{x_{3}}+x_{2}\partial
_{x_{4}}\right) F=0 \\
\widehat{X}_{6}^{\prime }F=\left( \lambda x_{1}+x_{2}\right) \partial
_{x_{1}}+\lambda x_{2}\partial _{x_{2}}+\left( \left( \lambda -\mu \right)
x_{3}+x_{4}\right) \partial _{x_{3}}+ \\
+\left( \lambda -\mu \right) x_{4}\partial _{x_{4}}=0
\end{array}
\right\}
\end{equation}
It is easy to see that $I_{1}=x_{2}\exp \left( \frac{-\lambda x_{1}}{x_{2}}%
\right) $ is a solution of the system, since it trivially verifies the first
equation. Consider the new variables $u=\ln \left(
x_{1}x_{4}-x_{2}x_{3}\right) $ and $v=x_{1}x_{2}^{-1}$. Then $\widehat{X}%
_{5}^{\prime }\left( u\right) =\widehat{X}_{5}^{\prime }\left( v\right) =0$
and $\widehat{X}_{6}^{\prime }u=\left( 2\lambda -\mu \right) ,\widehat{X}%
_{6}^{\prime }\left( v\right) =1$. Then the new equation
\begin{equation}
\partial _{u}F\left( u,v\right) +\frac{\widehat{X}_{6}^{\prime }\left(
v\right) }{\widehat{X}_{6}^{\prime }\left( u\right) }\partial _{v}F\left(
u,v\right) =0
\end{equation}
is known to provide a solution of $\left( 4\right) \;\cite{Dixon}$. Since $%
\left( 2\lambda -\mu \right) v-u$ is a solution of $\left( 5\right) $, the
system $\left( 4\right) $ has the solution $I_{2}=\left(
x_{1}x_{4}-x_{2}x_{3}\right) \exp \left( \frac{-\left( 2\lambda -\mu \right)
x_{1}}{x_{2}}\right) $. Clearly $I_{1}$ and $I_{2}$ are functionally
independent. and since $\mathcal{N}\left( \frak{g}_{6,65}\right) =2$, we
have found a fundamental set of solutions.
\end{example}

\section{Six dimensional Lie algebras}

The problem of obtaining criteria to classify solvable Lie algebras was
initiated by Mubarakzyanov in 1962 \cite{Mu1}, and succesfully applied to
the five dimensional case \cite{Mu2}. Solvable Lie algebras $\frak{r}$ admit
a decomposition $\frak{r}=\frak{n}\oplus \frak{t}$, where $\frak{n}$ is the
nilradical and $\left[ \frak{t},\frak{t}\right] \subset \frak{n},\;\left[
\frak{t},\frak{n}\right] \subset \frak{n}$. Since the nilradical must
satisfy $\dim \frak{n}\geq \frac{1}{2}\dim \frak{r}$ \cite{Mu1}, the
classification of six dimensional Lie algebras reduces to three cases,
namely, solvable Lie algebras having nilradical of dimension four, five or
three, the latter being immediate. The case $\dim \,\frak{n}=5$ was solved
by Mubarakzyanov in \cite{Mu3}, while the classification of the case $\dim \,%
\frak{n}=4$, which finishes the classification of six dimensional Lie
algebras, was presented by Turkowski in \cite{Tu1}$.$ The invariants of
solvable Lie algebras in dimension $n\leq 5$ and of nilpotent Lie algebras
in dimension $n\leq 6$ were computed in \cite{PWZ}, while the invariants of
six dimensional solvable algebras with four dimensional nilradical were
determined in \cite{Nd}. In this work, the author postulates that the
structure of invariants for six dimensional Lie algebras with a five
dimensional nilradical is rather similar to that of five dimensional
solvable algebras with four dimensional radical. This assertion is right in
what concerns the general structure of the invariants, but does not consider
the fact that even dimensional Lie algebras are allowed to have no
non-trivial invariant, as in fact happens. The analysis of these classes can
be generalized to obtain some criteria for the nonexistence of non-trivial
invariants in arbitrary dimension. This fact makes the analysis of these
invariants interesting, and finishes the determination of the invariants of
real or solvable Lie algebras in dimension $n\leq 6$.

The classification of the solvable algebras $\frak{r}$ such that $\dim \frac{%
\frak{r}}{\frak{n}}=1$ is based on the method developed in \cite{Mu4}, and
provides 99 isomorphism classes. These algebras are ordered and labelled
respect to the structure of its nilradical. The representatives of these
five dimensional nilpotent Lie algebras are taken from \cite{Mo}. In listing
the six dimensional Lie algebras of \cite{Mu3} and its invariants, we will
preserve the nomenclature used there for both the nilradicals and the
algebras. The following table presents the nonzero brackets for the nilradicals:

\begin{table}[h]
\caption{Nilradicals in dimension five}
\begin{tabular}
[c]{lllll}%
Name & Brackets & & & \\\hline
$5\frak{g}_{1}:$ & abelian algebra &  &  &  \\
$\frak{g}_{3,1}+2\frak{g}_{1}$ & $\left[ X_{2},X_{3}\right] =X_{1}.$ &  &  &
\\
$\frak{g}_{4,1}+\frak{g}_{1}$ & $\left[ X_{1},X_{5}\right] =X_{2},$ & $\left[
X_{4},X_{5}\right] =X_{1}.$ &  &  \\
$\frak{g}_{5,1}$ & $\left[ X_{3},X_{5}\right] =X_{1},$ & $\left[ X_{4},X_{5}%
\right] =X_{2}.$ &  &  \\
$\frak{g}_{5,2}$ & $\left[ X_{2},X_{5}\right] =X_{1},$ & $\left[ X_{3},X_{5}%
\right] =X_{2},$ & $\left[ X_{4},X_{5}\right] =X_{3}.$ &  \\
$\frak{g}_{5,3}$ & $\left[ X_{2},X_{4}\right] =X_{3},$ & $\left[ X_{2},X_{5}%
\right] =X_{1},$ & $\left[ X_{4},X_{5}\right] =X_{2}.$ &  \\
$\frak{g}_{5,4}$ & $\left[ X_{2},X_{4}\right] =X_{1},$ & $\left[ X_{3},X_{5}%
\right] =X_{1}.$ &  &  \\
$\frak{g}_{5,5}$ & $\left[ X_{3},X_{4}\right] =X_{1},$ & $\left[ X_{2},X_{5}%
\right] =X_{1},$ & $\left[ X_{3},X_{5}\right] =X_{2}.$ &  \\
$\frak{g}_{5,6}$ & $\left[ X_{3},X_{4}\right] =X_{1},$ & $\left[ X_{2},X_{5}%
\right] =X_{1},$ & $\left[ X_{3},X_{5}\right] =X_{2},$ & $\left[ X_{4},X_{5}%
\right] =X_{3}$\\\hline
\end{tabular}
\end{table}

\begin{remark}
The original list of Mubarakzyanov \cite{Mu3} presents some minor misprints, as well as the omission of some algebras. Here we present corrected lists, taking into account the
corrections indicated by Turkowski in \cite{Tu2} and the author.
\end{remark}

\subsection{Nilradical isomorphic to $5\frak{g}_{1}$}

There are 12 isomorphim classes of solvable Lie algebras having an abelian
nilradical, all of them depending on one, two or three parameters, with the
exception of $\frak{g}_{6,5}$. The invariants of solvable Lie algebras with
abelian nilradical are one of the few cases where a general theorem can be
established. Their number is determined by the following result:

\begin{proposition}
\cite{Nd2} If $\frak{r}$ is a solvable Lie algebra with abelian nilradical $%
\frak{n}$, then
\begin{equation*}
\mathcal{N}(\frak{r})=2\dim \,\frak{n}-\dim \,\frak{r}
\end{equation*}
\end{proposition}

Thus all the algebras in tables 1 and 2 admit four invariants. Moreover,
system $\left( 2\right) $ reduces considerably in this case, and the
invariants are gained from an unique equation. The fundamental set of
invariants are formed by rational functions whenever the element $X_{6}$
acts diagonally \cite{Nd}, and this happens only for $\frak{g}_{6,1}$. The
remaining algebras have at least one invariant of exponential form, although
some invariants reduce to rational or polynomial functions for specific
values of the parameters. In particular, if an element $X_{i}$ belongs to
the center of the Lie algebra, then the polynomial $x_{i}$ is an invariant,
since it is in the center of the universal enveloping algebra.

\subsection{Nilradical isomorphic to $\frak{g}_{3,1}+2\frak{g}_{1}$}

The case of a nilradical isomorphic to $\frak{g}_{3,1}+2\frak{g}_{1}$
provides the highest number of isomorphism classes, 25 in total [tables 4-6]. It is
immediately seen that the nilradical is isomorphic to the direct product of
the 3 dimensional Heisenberg Lie algebra with $\mathbb{K}^{2}$. This type of
maximal nilpotent ideal will tehrefore have some analogy to the other
general class of Lie algebras whose invariants are perfectly determined,
(the nilradical being isomorphic to the Heisenberg Lie algebra \cite{Wi})%
. In the present case, it is immediate from the tables that all Lie algebras
$\frak{g}_{6,13},..,\frak{g}_{6,38}$ have exactly 2 invariants, of which at
least one can be chosen as a rational function. An ocular inspection of
tables 4.-6. shows that the determination of the invariants for the algebras
having $\frak{g}_{3,1}+2\frak{g}_{1}$ as nilradical reduces to solve the
unique equation:

\begin{equation}
\left( C_{16}^{k}x_{k}\partial _{x_{1}}+C_{46}^{k}x_{k}
\partial _{x_{4}}+C_{56}^{k}x_{k}\partial _{x_{5}}\right) F=0
\end{equation}
which shows that the structure of the nilradical merely implies that an
invariant does not depend on the variables $x_{2}$ and $x_{3}$. From this we
will obtain generally a fundamental set of invariants formed by a rational
function and an exponential function.

\subsection{Nilradical isomorphic to $\frak{g}_{4,1}+\frak{g}_{1}$}

14 isomorphism classes have $\frak{g}_{4,1}+\frak{g}_{1}$ as nilradical, six of them depending on a continuous parameter [tables 7 and 8].
The nilradical $\frak{g}_{4,1}+\frak{g}_{1}\frak{\ }$is easily seen to be
the direct product of $\mathbb{K}$ and a four dimensional filiform Lie
algebra ( i.e., a nilpotent Lie algebra with maximal nilpotence index). The
system $\left( 2\right) $ reduces to the following system of two equations:
\begin{equation}
\left.
\begin{array}{r}
\left( x_{2}\partial _{x_{1}}+x_{1}\partial _{x_{4}}\right) F=0 \
\left( C_{16}^{k}x_{k} \partial _{x_{1}}+ C_{26}^{k}x_{k}
\partial _{x_{2}}+ C_{36}^{k}x_{k} \partial _{x_{3}}+ C_{46}^{k}x_{k} \partial _{x_{4}}\right) F=0
\end{array}
\right\}
\end{equation}
This tells that in general we will obtain a fundamental set of invariants
formed by an exponential and a rational function, depending on the action of
$X_{6}$ on the nilradical. For concrete values of the parameters the
exponential function can reduce to a rational, even to a polynomial
invariant.

We can deduce a general result starting from this case: let $\frak{g}_{2m}$
be the $2m$-dimensional nilpotent Lie algebra given by the brackets $\left[
X_{1},X_{i}\right] =X_{i+1},2\leq i\leq 2m-1$ over the basis $\left\{
X_{1},..,X_{2m}\right\} $. Consider a $\left( 2m+2\right) $-dimensional
solvable Lie algebra $\frak{r}_{2m+2}$\ \ whose nilradical is isomorphic to $%
\frak{g}_{2m}\oplus \mathbb{K}$. This generalizes in natural manner the Lie
algebras $\frak{g}_{6,53}\,$to $\frak{g}_{6,70}$.

\begin{theorem}
For $m\geq 2$ the number of invariants of $\frak{r}_{2m+2}$ is given by $%
\mathcal{N}\left( \frak{r}_{2m+2}\right) =2m-2$. Moreover, whenever $\left[
X_{2m},X_{2m+2}\right] \neq 0$ or $\left[ X_{2m+1},X_{2m+2}\right] \neq 0$,
an invariant $F$ is a solution of the following system:
\begin{equation}
\left.
\begin{array}{r}
\left( \sum_{j=2}^{2m-1}x_{j+1}\partial _{x_{j}}\right) F=0 \
\left( \sum_{i,j=2}^{2m+1}C_{i,2m+2}^{k}x_{k}\partial _{x_{i}}\right) F=0
\end{array}
\right\}
\end{equation}
\end{theorem}

\begin{proof}
The number of invariants follows at once from proposition 1, since the
matrix $\left( C_{ij}^{k}x_{k}\right) _{1\leq i<j\leq 2m+2}$ representing the
commutator table has rank four for any $m$. For the second assertion,
observe that since $X_{2m}$ and $X_{2m+1}$ belong to the center of the
nilradical, the associated differential operators are
\begin{equation}
\left.
\begin{array}{l}
\widehat{X}_{2m}=-C_{2m,2m+2}^{k}x_{k}\partial _{x_{2m+2}} \
\widehat{X}_{2m+1}=-C_{2m+1,2m+2}^{k}x_{k}\partial _{x_{2m+2}}
\end{array}
\right\}
\end{equation}
and if either   $\left[ X_{2m},X_{2m+2}\right] \neq 0$ or $\left[
X_{2m+1},X_{2m+2}\right] \neq 0$, this implies that $\partial _{x_{2m+2}}F=0$
for any invariant. On the other hand, from $\widehat{X}_{2m-1}=x_{2m}%
\partial _{x_{1}}-C_{2m-1,2m+2}^{k}x_{k}\partial _{x_{2m+2}}$ it also
follows that $\partial _{x_{1}}F=0$. This reduces system $\left( 2\right) $
to the previous form.
\end{proof}

The practical procedure to solve systems like $\left(5\right) $ is to solve
the first equation and then search for solutions which are also common to
the second equations \cite{Ol}.

\subsection{Nilradical isomorphic to $\frak{g}_{5,1}$}

There are 18 isomorphism classes of six dimensional Lie algebras with
nilradical isomorphic to $\frak{g}_{5,1}$ [tables 9 and 10]. From these, only the algebras $\frak{g}_{6,53}, \frak{g}_{6,66}$ and $\frak{g}_{6,69}$ do not depend on parameters. $\frak{g}_{6,53}$ is also the only algebra which admits a fundamental set of invariants formed by two polynomials. For the remaining cases, such a set is either formed by a rational and an exponential function or two exponential functions.

\subsection{Nilradical isomorphic to $\frak{g}_{5,2}$}

There are five isomorphism classes with $\frak{g}_{5,2}$ as nilradical, from which only one, $\frak{g}_{6,71}$, depends on a continuous parameter, while the algebra $\frak{g}_{6,73}$ depends on a parameter having the values $\pm 1$. It follows immediately from the brackets listed in table 11. that the invariants of these algebras do not depend on the variables $x_{5}$ and $x_{6}$, and that any of these algebras satisfies $\mathcal{N}=2$. With the exception of the algebras $\frak{g}_{6,73}$ and $\frak{g}_{6,75}$, the action of $X_{6}$ on the nilradical is diagonal, which implies that we can always find a fundamental set of invariants formed by rational functions.

\subsection{Nilradical isomorphic to $\frak{g}_{5,3}$}

There are 6 isomorphism classes having $\frak{g}_{5,3}$ as
nilradical, one of them depending on a continuous parameter. It follows immediately
from the structure of these algebras that there are exactly two functionally
independent invariants. Only for $\frak{g}_{76}$ and $\frak{g}_{6,80}$ the
action of $X_{6}$ is diagonal, thus only for these two algebras we obtain a fundamental set of invariants formed by
rational functions. It is immediately seen from table 12. that system $%
\left( 2\right) $ reduces to the following:
\begin{equation}
\left.
\begin{array}{c}
x_{3}\partial _{x_{4}}+x_{1}\partial _{x_{5}}=0 \
x_{3}\partial _{x_{2}}-x_{2}\partial _{x_{5}}=0 \
x_{1}\partial _{x_{2}}+x_{2}\partial _{x_{4}}=0 \
\sum_{i=1}^{5}C_{i6}^{k}x_{k}\partial _{x_{i}}=0
\end{array}
\right\} .
\end{equation}
From this we deduce that any invariant can be represented as a function $%
F\left( x_{1},x_{3},\frac{2(x_{3}x_{5}-x_{1}x_{4})+x_{2}^{2}}{x_{3}}\right) $%
, the new variables being a fundamental set of solutions for the three first
equations, since in any case $\frac{\partial F}{\partial x_{6}}=0$. Even
more is true: we can choose one of the fundamental invariants as a solution
of the equation
\begin{equation}
ax_{1}\partial _{x_{1}}+\left( bx_{1}+cx_{3}\right) \partial _{x_{3}}=0
\end{equation}
depending only on $x_{1}\,$and $x_{3}$, where $a,b,c$ are constants. Up to the case where $a=c$ and $b=1$%
, which provides and exponential solution, the remaining cases have the
function $f=x_{1}^{-c}\left( \left( a-c\right) x_{3}-bx_{1}\right) ^{a}$ as
solution. This pattern occurs for all the algebras in this case, up to $\frak{%
g}_{6,79}$.

\subsection{Nilradical isomorphic to $\frak{g}_{5,4}$}

There are 12 isomorphism classes having this algebra as nilradical,
eight of them depending on one or more parameters [tables 13 and 14]. This is also the first
case where the number of invariants depends essentially on the action of the
vector $X_{6}$ on the nilradical $\frak{g}_{5,4}$. As follows from \cite{Mu3}%
, the structure of solvable Lie algebras having this algebra as nilradical
is
\begin{eqnarray*}
\left[ X_{2},X_{4}\right]  &=&X_{1};\;\left[ X_{3},X_{5}\right] =X_{1} \\
\left[ X_{i},X_{6}\right]  &=&C_{i6}^{k}x_{k},\;1\leq i,k\leq 5
\end{eqnarray*}
Therefore the matrix $\left( C_{ij}^{k}x_{k}\right) $ representing the commutator
table is
\begin{equation}
A=\left(
\begin{array}{cccccc}
0 & 0 & 0 & 0 & 0 & C_{16}^{k}x_{k} \\
0 & 0 & 0 & x_{1} & 0 & C_{26}^{k}x_{k} \\
0 & 0 & 0 & 0 & x_{1} & C_{36}^{k}x_{k} \\
0 & -x_{1} & 0 & 0 & 0 & C_{46}^{k}x_{k} \\
0 & 0 & -x_{1} & 0 & 0 & C_{56}^{k}x_{k} \\
-C_{16}^{k}x_{k} & -C_{26}^{k}x_{k} & -C_{36}^{k}x_{k} & -C_{46}^{k}x_{k} &
-C_{56}^{k}x_{k} & 0
\end{array}
\right)
\end{equation}
whose determinant is $d=x_{1}^{4}\left( C_{16}^{k}x_{k}\right) ^{2}$.
Therefore, we will obtain non-trivial invariants if and only if $C_{16}^{k}=0
$ for $1\leq k\leq 5$. If this occurs, $A$ is easily seen to have rank four
and by the Beltrametti-Blasi formula \cite{Be}\ the algebra will have
exactly two invariants. Moreover, since this implies that $X_{1}$ is a
central element, one invariant can be chosen as $I_{1}=x_{1}$.

\begin{lemma}

Let $\frak{g}$ be a solvable Lie algebra whose nilradical is isomorphic to $%
\frak{g}_{5,4}$. If $\det \left( A\right) =0$, then $\frak{g}$ admits a

fundamental set of invariants formed by two polynomial functions, $%
I_{1}=x_{1}$ being one of them.

\end{lemma}

Thus in this case we obtain classical Casimir operators \cite{AA,Di}. Observe further that $\frak{g}_{5,4}$ is isomorphic to the five dimensional Heisenberg Lie algebra. In \cite{Wi} the invariants of solvable algebras having a Heisenberg Lie algebra as nilradical are determined generically.

\subsection{Nilradical isomorphic to $\frak{g}_{5,5}$}

\bigskip Five isomorphim classes arise from Mubarakzyanovs classification
\cite{Mu3} for this nilradical, only one of them, $\frak{g}_{6,94}$,
depending on a parameter [table 15]. An argumentation like in the preceding case shows
that these algebras have no non-trivial invariants, up to the algebra $\frak{%
g}_{6,94}$ for the value $\lambda +2=0$, for which we obtain again a
fundamental set of invariants formed by two polynomial functions.

\subsection{Nilradical isomorphic to $\frak{g}_{5,6}$}

Finally, for the nilradical $\frak{g}_{5,6}$ there is a unique
non-parametrized isomorphism class, $\frak{g}_{6,99}$ [table 16]. For this algebra is
is easily verified that the matrix $\left( C_{ij}^{k}x_{k}\right) $ has maximal
rank, which shows that this algebra has no invariants.

The last three nilradicals, $\frak{g}_{5,4},\frak{g}_{5,5}$ and $\frak{g}%
_{5,6}$ have in common that they are all deformations of the 5-dimensional Heisenberg Lie algebra (considering the trivial deformation for the Lie algebra $\frak{g}_{5,4}$) \cite{Ca2}. This fact allows us to deduce a similar result which holds in arbitrary dimension: For
$m\leq 2$ let $\frak{n}_{\alpha _{1},..,\alpha _{\left[ \frac{m}{2}\right]
+2}}$ be the $\left( 2m+1\right) $-dimensional nilpotent Lie algebra given
by the brackets
\begin{equation}
\left.
\begin{tabular}{ll}
$\left[ X_{2+j},X_{2m+1-j}\right] =X_{1},$ & $0\leq j\leq m-1$ \\
$\left[ X_{2+j},X_{2m+1}\right] =\alpha _{j}X_{1+j},$ & $1\leq j\leq 2m-2$
\\
$\alpha _{j}+\alpha _{2m-j}=0,$ & $1\leq j\leq m-1$%
\end{tabular}
\right\} ,
\end{equation}
over the basis $\left\{ X_{1},..,X_{2m+1}\right\} $. Let $\frak{r}_{m}$ be a
$\left( 2m+2\right) $-dimensional solvable Lie algebra whose nilradical is
isomorphic to $\frak{n}_{\alpha _{1},..,\alpha _{\left[ \frac{m}{2}\right]
+2}}$.

\begin{theorem}
Let $\frak{r}_{m}$ be a $\left( 2m+2\right) $-dimensional solvable Lie
algebra whose nilradical is isomorphic to $\frak{n}_{\alpha _{1},..,\alpha _{%
\left[ \frac{m}{2}\right] +2}}$.

\begin{enumerate}

\item  If the action $\left( \frac{r_{m}}{\frak{n}_{\alpha _{1},..,\alpha _{%
\left[ \frac{m}{2}\right] +2}}}\right) $ over $X_{1}$ is nonzero, then $%
\mathcal{N}\left( \frak{r}_{m}\right) =0$.

\item  If $\left( \frac{r_{m}}{\frak{n}_{\alpha _{1},..,\alpha _{\left[
\frac{m}{2}\right] +2}}}\right) .X_{1}=0$, then $\frak{r}_{m}$ admits a
fundamental set of invariants formed by two elements, one of which can be
taken as $I_{1}=x_{1}$.
\end{enumerate}
\end{theorem}

\begin{proof}
If $\left( C_{ij}^{k}x_{k}\right) $ is the matrix associated to the commutators
of $\frak{r}_{m}$, then it can easily be proven by induction over $m$ that
\begin{equation*}
\det \left( C_{ij}^{k}\right) =x_{1}^{2m}\left( C_{1,2m+2}^{k}x_{k}\right)
^{2}
\end{equation*}
so that $\mathcal{N}\left( \frak{r}_{m}\right) =0$ if and only if $%
C_{1,2m+2}^{k}=0$ for all $k$. Otherwise we obtain that $rank\left(
C_{ij}^{k}\right) =2m$, thus the algebra has two invariants. The second
assertion follows from the fact that $X_{1}$ spans the center of the
nilradical.
\end{proof}

Although there is no procedure to determine the detailed structure of the
second invariant, the analysis of low values of $m$ seems to indicate that
it can always be chosen as a polynomial function. The exact structure of this polynomial depends heavily on the parameters $\alpha_{i}$, and cannot therefore be comprised in a formula.

\begin{table}[h]
\caption{Nilradical isomorphic to $5\mathfrak{g}_{1}$}
\begin{tabular}
[c]{llll}%
Name & Brackets & & Invariants\\\hline
$\frak{g}_{6,1}$ & $\left[  X_{1},X_{6}\right]  =X_{1}$ & $\left[  X_{2}%
,X_{6}\right]  =\alpha X_{2}$ & $I_{1}=x_{2}x_{1}^{-\alpha}$\\
$0<\left|  \delta\right|  \leq\left|  \gamma\right|$   & $\left[  X_{3},X_{6}\right]  =\beta X_{3}$ &
$\left[  X_{4},X_{6}\right]  =\gamma X_{4}$ & $I_{2}=x_{3}x_{1}^{-\beta}$\\
$\leq\left|  \beta\right|
\leq\left|  \alpha\right|  <1$& $\left[  X_{5},X_{6}\right]  =\delta X_{5}$ &  & $I_{3}=x_{4}x_{1}^{-\gamma
}$\\
&  &  & $I_{4}=x_{5}x_{1}^{-\delta}$\\
$\frak{g}_{6,2}$ & $\left[  X_{1},X_{6}\right]  =\alpha X_{1}$ & $\left[
X_{2},X_{6}\right]  =X_{1}+\alpha X_{2}$ & $I_{1}=x_{1}\exp\left(
\frac{-\alpha x_{2}}{x_{1}}\right)  $\\
$0<\left|  \delta\right|  \leq\left|  \gamma\right|  \leq1$ & $\left[
X_{3},X_{6}\right]  =X_{3}$ & $\left[  X_{4},X_{6}\right]  =\gamma X_{4}$ &
$I_{2}=x_{3}\exp\left(  \frac{-x_{2}}{x_{1}}\right)  $\\
& $\left[  X_{5},X_{6}\right]  =\delta X_{5}$ &  & $I_{3}=x_{4}\exp\left(
\frac{-\gamma x_{2}}{x_{1}}\right)  $\\
\hline

\end{tabular}%

\end{table}

\begin{table}[h]

\caption{Nilradical isomorphic to $5\mathfrak{g}_{1}$ (cont.)}

\begin{tabular}
[c]{llll}%
Name & Brackets & & Invariants\\\hline
$\frak{g}_{6,9}$ & $\left[  X_{1},X_{6}\right]  =\alpha X_{1}$ & $\left[
X_{2},X_{6}\right]  =\beta X_{2}$ & $I_{1}=\left(  x_{4}^{2}+x_{5}^{2}\right)
\left(  \frac{x_{4}+ix_{5}}{x_{4}-ix_{5}}\right)  ^{ip}$\\
$\alpha\neq0$ & $\left[  X_{3},X_{6}\right]  =X_{2}+\beta X_{3}$ & $\left[
X_{4},X_{6}\right]  =pX_{4}-X_{5}$ & $I_{2}=x_{1}\exp^{\alpha\arctan
x_{4}x_{5}^{-1}}$\\
& $\left[  X_{5},X_{6}\right]  =X_{4}+pX_{5}$ &  & $I_{3}=x_{2}\exp
^{\beta\arctan x_{4}x_{5}^{-1}}$\\
&  &  & $I_{4}=(x_{3}+x_{2}\arctan\frac{x_{4}}{x_{5}})\exp^{\beta\arctan\frac{
x_{4}}{x_{5}}}$\\

$\frak{g}_{6,10}$ & $\left[  X_{1},X_{6}\right]  =\alpha X_{1}$ & $\left[
X_{2},X_{6}\right]  =X_{1}+\alpha X_{2}$ & $I_{1}=\left(  x_{4}^{2}+x_{5}%
^{2}\right)  \left(  \frac{x_{4}+ix_{5}}{x_{4}-ix_{5}}\right)  ^{ip}$\\
& $\left[  X_{3},X_{6}\right]  =X_{2}+\alpha X_{3}$ & $\left[  X_{4}%
,X_{6}\right]  =pX_{4}-X_{5}$ & $I_{2}=x_{1}\exp^{\alpha\arctan x_{4}%
x_{5}^{-1}}$\\
& $\left[  X_{5},X_{6}\right]  =X_{4}+pX_{5}$ &  & $I_{3}=(x_{2}+x_{1}%
\arctan\frac{x_{4}}{x_{5}})\exp^{\alpha\arctan\frac{
x_{4}}{x_{5}}}$\\
\hline
\end{tabular}
\end{table}

\begin{table}[h]

\caption{Nilradical isomorphic to $\mathfrak{g}_{3,1}+2\mathfrak{g}_{1}$ }
\begin{tabular}
[c]{llll}%
Name & Brackets & & Invariants\\\hline

$\frak{g}_{6,13}$ & $\left[  X_{2},X_{3}\right]  =X_{1}$ & $\left[

X_{1},X_{6}\right]  =\left(  \alpha+\beta\right)  X_{1}$ & $I_{1}%

=x_{4}^{\left(  \alpha+\beta\right)  }x_{1}^{-1}$\\

$\alpha\neq0$ & $\left[  X_{2},X_{6}\right]  =\alpha X_{2}$ & $\left[

X_{3},X_{6}\right]  =\beta X_{3}$ & $I_{2}=x_{4}^{h}x_{5}^{-1}$\\

& $\left[  X_{4},X_{6}\right]  =X_{4}$ & $\left[  X_{5},X_{6}\right]  =hX_{5}$%

& \\

$\frak{g}_{6,14}$ & $\left[  X_{2},X_{3}\right]  =X_{1}$ & $\left[

X_{1},X_{6}\right]  =\left(  \alpha+\beta\right)  X_{1}$ & $I_{1}%

=x_{4}^{\left(  \alpha+\beta\right)  }x_{1}^{-1}$\\

& $\left[  X_{2},X_{6}\right]  =\alpha X_{2}$ & $\left[  X_{3},X_{6}\right]

=\beta X_{3}$ & $I_{2}=x_{1}\exp\left(  \frac{-\left(  \alpha+\beta\right)

x_{5}}{x_{1}}\right)  $\\

& $\left[  X_{4},X_{6}\right]  =X_{4}$ & $\left[  X_{5},X_{6}\right]

=X_{1}+\left(  \alpha+\beta\right)  X_{5}$ & \\

$\frak{g}_{6,15}$ & $\left[  X_{2},X_{3}\right]  =X_{1}$ & $\left[

X_{1},X_{6}\right]  =\left(  1+h\right)  X_{1}$ & $I_{1}=x_{4}^{\left(

1+h\right)  }x_{1}^{-1}$\\
\hline

\end{tabular}

\end{table}

\begin{table}[h]

\caption{Nilradical isomorphic to $\mathfrak{g}_{3,1}+2\mathfrak{g}_{1}$ (cont.)}

\begin{tabular}
[c]{llll}%
Name & Brackets & & Invariants\\\hline

$\frak{g}_{6,24}$ & $\left[  X_{2},X_{3}\right]  =X_{1}$ & $\left[

X_{1},X_{6}\right]  =0$ & $I_{1}=x_{1}$\\

& $\left[  X_{2},X_{6}\right]  =X_{3}$ & $\left[  X_{3},X_{6}\right]  =X_{4}$%

& $I_{2}=x_{5}^{h}\exp\left(  \frac{-x_{4}}{x_{1}}\right)  $\\

& $\left[  X_{4},X_{6}\right]  =hX_{1}$ & $\left[  X_{5},X_{6}\right]  =X_{5}$%

& \\

$\frak{g}_{6,25}$ & $\left[  X_{2},X_{3}\right]  =X_{1}$ & $\left[

X_{1},X_{6}\right]  =\left(  1+h\right)  X_{1}$ & $I_{1}=x_{1}^{-\beta}%

x_{5}^{1+h}$\\

& $\left[  X_{2},X_{6}\right]  =X_{2}$ & $\left[  X_{3},X_{6}\right]  =hX_{3}$%

& $I_{2}=x_{1}\exp\left(  \frac{-\left(  1+h\right)  x_{4}}{x_{5}}\right)  $\\

& $\left[  X_{4},X_{6}\right]  =\beta X_{4}+X_{5}$ & $\left[  X_{5}%

,X_{6}\right]  =\beta X_{5}$ & \\

$\frak{g}_{6,26}$ & $\left[  X_{2},X_{3}\right]  =X_{1}$ & $\left[

X_{1},X_{6}\right]  =\left(  1+h\right)  X_{1}$ & $I_{1}=\frac{2x_{1}%

x_{4}-x_{5}^{2}}{x_{1}^{2}}$\\
\hline

\end{tabular}

\end{table}

\begin{table}[h]

\caption{Nilradical isomorphic to $\mathfrak{g}_{3,1}+2\mathfrak{g}_{1}$ (cont.)}

\begin{tabular}
[c]{llll}%
Name & Brackets & & Invariants\\\hline

$\frak{g}_{6,34}$ & $\left[  X_{2},X_{3}\right]  =X_{1}$ & $\left[

X_{1},X_{6}\right]  =2\alpha X_{1}$ & $I_{1}=$\\

$\varepsilon h=0$ & $\left[  X_{2},X_{6}\right]  =\alpha X_{2}+X_{3}$ &

$\left[  X_{3},X_{6}\right]  =-X_{2}+\alpha X_{3}$ & $I_{2}=$\\

& $\left[  X_{4},X_{6}\right]  =\left(  2\alpha+h\right)  X_{4}+X_{5}$ &

$\left[  X_{5},X_{6}\right]  =\left(  2\alpha+h\right)  X_{5}+\varepsilon

X_{1}$ & \\

$\frak{g}_{6,35}$ & $\left[  X_{2},X_{3}\right]  =X_{1}$ & $\left[

X_{1},X_{6}\right]  =\left(  \alpha+\beta\right)  X_{1}$ & $I_{1}=\left(

x_{4}^{2}+x_{5}^{2}\right)  \left(  \frac{x_{4}-ix_{5}}{x_{4}+ix_{5}}\right)

^{i\gamma}$\\

$_{\alpha^{2}+\beta^{2}\neq 0}$ & $\left[  X_{2},X_{6}\right]  =\alpha X_{2}$ &

$\left[  X_{3},X_{6}\right]  =\beta X_{3}$ & $I_{2}=x_{1}\exp\left(  -\left(

\alpha+\beta\right)  \arctan\left(  \frac{x_{4}}{x_{5}}\right)  \right)  $\\

& $\left[  X_{4},X_{6}\right]  =\gamma X_{4}+X_{5}$ & $\left[  X_{5}%

,X_{6}\right]  =-X_{4}+\gamma X_{5}$ & \\

$\frak{g}_{6,36}$ & $\left[  X_{2},X_{3}\right]  =X_{1}$ & $\left[

X_{1},X_{6}\right]  =2\alpha X_{1}$ & $I_{1}=\left(  x_{4}^{2}+x_{5}%

^{2}\right)  \left(  \frac{x_{4}-ix_{5}}{x_{4}+ix_{5}}\right)  ^{i\gamma}$\\
\hline

\end{tabular}

\end{table}

\begin{table}[h]

\caption{Nilradical isomorphic to $\mathfrak{g}_{4,1}+\mathfrak{g}_{1}$}

\begin{tabular}
[c]{llll}%

Name & Brackets & & Invariants\\\hline

$\frak{g}_{6,39}$ & $\left[  X_{1},X_{5}\right]  =X_{2}$ & $\left[

X_{4},X_{5}\right]  =X_{1}$ & $I_{1}=x_{3}^{2+h}x_{2}^{-\gamma}$\\

$\gamma\neq0$ & $\left[  X_{1},X_{6}\right]  =\left(  1+h\right)  X_{1}$ &

$\left[  X_{2},X_{6}\right]  =\left(  2+h\right)  X_{2}$ & $I_{2}=\frac

{x_{1}^{2}-2x_{2}x_{4}}{x_{2}x_{3}\frac{h}{\gamma}}$\\

& $\left[  X_{3},X_{6}\right]  =\gamma X_{3}$ & $\left[  X_{4},X_{6}\right]

=hX_{4}$ & \\

& $\left[  X_{5},X_{6}\right]  =X_{5}$ &  & \\

$\frak{g}_{6,40}$ & $\left[  X_{1},X_{5}\right]  =X_{2}$ & $\left[

X_{4},X_{5}\right]  =X_{1}$ & $I_{1}=x_{2}\exp\left(  \frac{-\left(

2+h\right)  x_{3}}{x_{2}}\right)  $\\

& $\left[  X_{1},X_{6}\right]  =\left(  1+h\right)  X_{1}$ & $\left[

X_{2},X_{6}\right]  =\left(  2+h\right)  X_{2}$ & $I_{2}=\frac{x_{1}%

^{2}-2x_{2}x_{4}}{x_{2}}\exp\left(  \frac{-hx_{3}}{x_{2}}\right)  $\\

& $\left[  X_{3},X_{6}\right]  =X_{2}+\left(  2+h\right)  X_{3}$ & $\left[

X_{4},X_{6}\right]  =hX_{4}$ & \\
\hline

\end{tabular}

\end{table}

\begin{table}[h]

\caption{Nilradical isomorphic to $\mathfrak{g}_{3,1}+2\mathfrak{g}_{1}$ (cont.)}

\begin{tabular}
[c]{llll}%

Name & Brackets & & Invariants\\\hline

$\frak{g}_{6,49}$ & $\left[  X_{1},X_{5}\right]  =X_{2}$ & $\left[

X_{4},X_{5}\right]  =X_{1}$ & $I_{1}=x_{3}$\\

$\epsilon\neq 0, \pm1$& $\left[  X_{1},X_{6}\right]  =X_{1}$ & $\left[  X_{2},X_{6}\right]  =X_{2}$%

& $I_{2}=x_{2}^{\varepsilon}\exp\left(  \frac{-2x_{2}x_{4}+x_{1}^{2}}%

{2x_{2}^{2}}\right)  $\\

& $\left[  X_{3},X_{6}\right]  =0$ & $\left[  X_{4},X_{6}\right]  =\varepsilon

X_{2}+X_{4}$ & \\

& $\left[  X_{5},X_{6}\right]  =X_{3}$ &  & \\

$\frak{g}_{6,50}$ & $\left[  X_{1},X_{5}\right]  =X_{2}$ & $\left[

X_{4},X_{5}\right]  =X_{1}$ & $I_{1}=x_{2}^{\varepsilon}\exp\left(

\frac{-x_{3}}{x_{2}}\right)  $\\

$\epsilon\neq 0, \pm1$& $\left[  X_{1},X_{6}\right]  =X_{1}$ & $\left[  X_{2},X_{6}\right]  =X_{2}$%

& $I_{2}=x_{2}^{-2}\varepsilon\left(-2x_{2}x_{4}\right)+x_{3}^{2} $\\

& $\left[  X_{3},X_{6}\right]  =\varepsilon X_{2}+X_{3}$ & $\left[

X_{4},X_{6}\right]  =X_{3}+X_{4}$ & \\

& $\left[  X_{5},X_{6}\right]  =0$ &  & \\

$\frak{g}_{6,51}$ & $\left[  X_{1},X_{5}\right]  =X_{2}$ & $\left[

X_{4},X_{5}\right]  =X_{1}$ & $I_{1}=x_{2}$\\

$\epsilon\neq \pm1$& $\left[  X_{1},X_{6}\right]  =0$ & $\left[  X_{2},X_{6}\right]  =0$ &

$I_{2}=x_{3}^{\varepsilon}\exp\left(  \frac{2x_{2}x_{4}-x_{1}^{2}}{x_{2}%

}\right)  $\\

& $\left[  X_{3},X_{6}\right]  =X_{3}$ & $\left[  X_{4},X_{6}\right]

=\varepsilon X_{2}$ & \\

& $\left[  X_{5},X_{6}\right]  =0$ &  & \\

$\frak{g}_{6,52}$ & $\left[  X_{1},X_{5}\right]  =X_{2}$ & $\left[

X_{4},X_{5}\right]  =X_{1}$ & $I_{1}=x_{2}$\\

$\epsilon\neq 0, \pm1$& $\left[  X_{1},X_{6}\right]  =0$ & $\left[  X_{2},X_{6}\right]  =0$ &

$I_{2}=x_{3}^{\varepsilon}\exp\left(  \frac{2x_{2}x_{4}-x_{1}^{2}}{x_{2}%

}\right)  $\\

& $\left[  X_{3},X_{6}\right]  =X_{3}$ & $\left[  X_{4},X_{6}\right]

=\varepsilon X_{2}$ & \\

& $\left[  X_{5},X_{6}\right]  =X_{4}$ &  & \\\hline

\end{tabular}

\end{table}

\begin{table}[h]

\caption{Nilradical isomorphic to $\mathfrak{g}_{5,1}$}
\begin{tabular}
[c]{llll}%
Name & Brackets & & Invariants\\\hline

$\frak{g}_{6,53}$ & $\left[  X_{3},X_{5}\right]  =X_{1}$ & $\left[

X_{4},X_{5}\right]  =X_{2}$ & $I_{1}=x_{1}$\\

& $\left[  X_{1},X_{6}\right]  =0$ & $\left[  X_{2},X_{6}\right]  =0$ &

$I_{2}=x_{2}$\\

& $\left[  X_{3},X_{6}\right]  =X_{3}$ & $\left[  X_{4},X_{6}\right]  =X_{4}$%

& \\

& $\left[  X_{5},X_{6}\right]  =-X_{5}$ &  & \\

$\frak{g}_{6,54}$ & $\left[  X_{3},X_{5}\right]  =X_{1}$ & $\left[

X_{4},X_{5}\right]  =X_{2}$ & $I_{1}=x_{2}x_{1}^{-\lambda}$\\

& $\left[  X_{1},X_{6}\right]  =X_{1}$ & $\left[  X_{2},X_{6}\right]  =\lambda

X_{2}$ & $I_{2}=\left(  x_{1}x_{4}-x_{2}x_{3}\right)  x_{1}^{\alpha+1-\gamma}%

$\\

& $\left[  X_{3},X_{6}\right]  =\left(  1-\gamma\right)  X_{3}$ & $\left[

X_{4},X_{6}\right]  =\left(  \lambda-\gamma\right)  X_{4}$ & \\
\hline

\end{tabular}

\end{table}

\begin{table}[h]

\caption{Nilradical isomorphic to $\mathfrak{g}_{5,1}$ (cont.)}

\begin{tabular}
[c]{llll}%

Name & Brackets & & Invariants\\\hline

$\frak{g}_{6,62}$ & $\left[  X_{3},X_{5}\right]  =X_{1}$ & $\left[

X_{4},X_{5}\right]  =X_{2}$ & $I_{1}=x_{2}^{2}x_{1}^{-1}$\\

$\omega=0,1$ & $\left[  X_{1},X_{6}\right]  =2X_{1}$ & $\left[  X_{2}%

,X_{6}\right]  =2X_{2}$ & $I_{2}=x_{1}\exp\left(  \frac{2\left(  x_{1}%

x_{4}-x_{2}x_{3}\right)  }{x_{2}^{2}}\right)  $\\

& $\left[  X_{3},X_{6}\right]  =X_{2}+X_{3}$ & $\left[  X_{4},X_{6}\right]

=0$ & \\

& $\left[  X_{5},X_{6}\right]  =\omega X_{3}+X_{5}$ &  & \\

$\frak{g}_{6,63}$ & $\left[  X_{3},X_{5}\right]  =X_{1}$ & $\left[

X_{4},X_{5}\right]  =X_{2}$ & $I_{1}=x_{2}x_{1}^{-\lambda}$\\

& $\left[  X_{1},X_{6}\right]  =X_{1}$ & $\left[  X_{2},X_{6}\right]  =\lambda

X_{2}$ & $I_{2}=x_{1}\exp\left(  \frac{-x_{1}x_{4}+x_{2}x_{3}}{x_{1}x_{2}%

}\right)  $\\

& $\left[  X_{3},X_{6}\right]  =X_{3}$ & $\left[  X_{4},X_{6}\right]

=X_{2}+\lambda X_{4}$ & \\
\hline

\end{tabular}

\end{table}

\begin{table}[h]

\caption{Nilradical isomorphic to $\mathfrak{g}_{5,2}$ }

\begin{tabular}
[c]{llll}%

Name & Brackets & & Invariants\\\hline

$\frak{g}_{6,71}$ & $\left[  X_{2},X_{5}\right]  =X_{1}$ & $\left[

X_{3},X_{5}\right]  =X_{2}$ & $I_{1}=\frac{\left(  3x_{1}^{2}x_{4}+x_{2}%

^{3}-3x_{1}x_{2}x_{3}\right)  ^{h+1}}{x_{1}^{h+2}\left(  2x_{1}x_{3}-x_{2}%

^{2}\right)  ^{h}}$\\

& $\left[  X_{4},X_{5}\right]  =X_{3}$ & $\left[  X_{1},X_{6}\right]  =\left(

h+3\right)  X_{1}$ & $I_{2}=\left(  -2x_{3}+x_{2}^{2}x_{1}^{-1}\right)

^{h+3}x_{1}^{-\left(  h+1\right)  }$\\

& $\left[  X_{2},X_{6}\right]  =\left(  h+2\right)  X_{2}$ & $\left[

X_{3},X_{6}\right]  =\left(  h+1\right)  X_{3}$ & \\

& $\left[  X_{4},X_{6}\right]  =hX_{4}$ & $\left[  X_{5},X_{6}\right]  =X_{5}$%

& \\

$\frak{g}_{6,72}$ & $\left[  X_{2},X_{5}\right]  =X_{1}$ & $\left[

X_{3},X_{5}\right]  =X_{2}$ & $I_{1}=\frac{\left(  3x_{1}^{2}x_{4}+x_{2}%

^{3}-3x_{1}x_{2}x_{3}+x_{2}^{3}\right)  ^{4}}{x_{1}^{9}}$\\

& $\left[  X_{4},X_{5}\right]  =X_{3}$ & $\left[  X_{1},X_{6}\right]  =4X_{1}$%

& $I_{2}=\left(  -2x_{1}x_{3}+x_{2}^{2}\right)  ^{2}x_{1}^{-3}$\\

& $\left[  X_{2},X_{6}\right]  =3X_{2}$ & $\left[  X_{3},X_{6}\right]

=2X_{3}$ & \\
\hline

\end{tabular}

\end{table}

\begin{table}[h]

\caption{Nilradical isomorphic to $\mathfrak{g}_{5,3}$ }

\begin{tabular}
[c]{llll}%

Name & Brackets & & Invariants\\\hline

$\frak{g}_{6,76}$ & $\left[  X_{2},X_{4}\right]  =X_{3}$ & $\left[

X_{2},X_{5}\right]  =X_{1}$ & $I_{1}=x_{3}^{\left(  2h+1\right)  }%

x_{1}^{-\left(  h+2\right)  }$\\

& $\left[  X_{4},X_{5}\right]  =X_{2}$ & $\left[  X_{1},X_{6}\right]  =\left(

2h+1\right)  X_{1}$ & $I_{2}=\frac{2x_{3}x_{5}+x_{2}^{2}-2x_{4}x_{1}^{\frac

{h}{2h+1}}}{x_{3}x_{1}^{\frac{h}{2h+1}}}$\\

& $\left[  X_{2},X_{6}\right]  =\left(  h+1\right)  X_{2}$ & $\left[

X_{3},X_{6}\right]  =\left(  h+2\right)  X_{3}$ & \\

& $\left[  X_{4},X_{6}\right]  =X_{4}$ & $\left[  X_{5},X_{6}\right]  =hX_{5}$%

& \\

$\frak{g}_{6,77}$ & $\left[  X_{2},X_{4}\right]  =X_{3}$ & $\left[

X_{2},X_{5}\right]  =X_{1}$ & $I_{1}=x_{1}^{\varepsilon}\exp\left(

\frac{-2x_{3}x_{5}-x_{2}^{2}+2x_{1}x_{4}}{2x_{1}^{2}}\right)  $\\

$\varepsilon=\pm1$ & $\left[  X_{4},X_{5}\right]  =X_{2}$ & $\left[

X_{1},X_{6}\right]  =X_{1}$ & $I_{2}=x_{3}x_{1}^{-2}$\\

& $\left[  X_{2},X_{6}\right]  =X_{2}$ & $\left[  X_{3},X_{6}\right]  =2X_{3}$%

& \\
\hline

\end{tabular}

\end{table}

\begin{table}[h]

\caption{Nilradical isomorphic to $\mathfrak{g}_{5,4}$ }

\begin{tabular}
[c]{llll}%

Name & Brackets & & Invariants\\\hline

$\frak{g}_{6,82}  $ & $\left[  X_{2},X_{4}\right]

=X_{1}$ & $\left[  X_{3},X_{5}\right]  =X_{1}$ & none\\

$\alpha=2$& $\left[  X_{1},X_{6}\right]  =2X_{1}$ & $\left[  X_{2},X_{6}\right]

=\left(  \lambda+1\right)  X_{2}$ & \\

& $\left[  X_{3},X_{6}\right]  =\left(  \lambda_{1}+1\right)  X_{3}$ &

$\left[  X_{4},X_{6}\right]  =\left(  -\lambda+1\right)  X_{4}$ & \\

& $\left[  X_{5},X_{6}\right]  =\left(  -\lambda_{1}+1\right)  X_{5}$ &  & \\

$\frak{g}_{6,82}$ & $\left[  X_{2},X_{4}\right]

=X_{1}$ & $\left[  X_{3},X_{5}\right]  =X_{1}$ & $I_{1}=x_{1}$\\

$\alpha=0$& $\left[  X_{1},X_{6}\right]  =X_{1}$ & $\left[  X_{2},X_{6}\right]  =\lambda

X_{2}$ & $I_{2}=\lambda x_{2}x_{4}+\lambda_{1}x_{3}x_{5}$\\

& $\left[  X_{3},X_{6}\right]  =\lambda_{1}X_{3}$ & $\left[  X_{4}%

,X_{6}\right]  =-\lambda X_{4}$ & $-x_{1}x_{6}$\\
\hline

\end{tabular}

\end{table}

\begin{table}[h]

\caption{Nilradical isomorphic to $\mathfrak{g}_{5,4}$ (cont.)}

\begin{tabular}
[c]{llll}%

Name & Brackets & & Invariants\\\hline

$\frak{g}_{6,89} $ & $\left[  X_{2},X_{4}\right]

=X_{1}$ & $\left[  X_{3},X_{5}\right]  =X_{1}$ & none\\

$\alpha\neq0$& $\left[  X_{1},X_{6}\right]  =\alpha X_{1}$ & $\left[  X_{2},X_{6}\right]

=\left(  \frac{\alpha}{2}+s\right)  X_{2}$ & \\

& $\left[  X_{3},X_{6}\right]  =\upsilon_{0}X_{5}+\frac{\alpha}{2}X_{3}$ &

$\left[  X_{4},X_{6}\right]  =\left(  \frac{\alpha}{2}-s\right)  X_{4}$ & \\

& $\left[  X_{5},X_{6}\right]  =-\upsilon_{0}X_{3}+\frac{\alpha}{2}X_{5}$ &

& \\

$\frak{g}_{6,89}$ & $\left[  X_{2},X_{4}\right]

=X_{1}$ & $\left[  X_{3},X_{5}\right]  =X_{1}$ & $I_{1}=x_{1}$\\

$\alpha=0$& $\left[  X_{1},X_{6}\right]  =0$ & $\left[  X_{2},X_{6}\right]  =sX_{2}$ &

$I_{2}=\frac{2x_{1}x_{6}-2sx_{2}x_{4}-\upsilon_{0}\left(  x_{3}^{2}+x_{5}%

^{2}\right)}{2}  $\\

& $\left[  X_{3},X_{6}\right]  =\upsilon_{0}X_{5}$ & $\left[  X_{4}%

,X_{6}\right]  =-sX_{4}$ & \\
\hline

\end{tabular}

\end{table}

\begin{table}[h]

\caption{Nilradical isomorphic to $\mathfrak{g}_{5,5}$ }

\begin{tabular}
[c]{llll}%

Name & Brackets & & Invariants\\\hline

$\frak{g}_{6,94} $ & $\left[  X_{3}%

,X_{4}\right]  =X_{1}$ & $\left[  X_{2},X_{5}\right]  =X_{1}$ & none\\

$\lambda+2\neq0$& $\left[  X_{3},X_{5}\right]  =X_{2}$ & $\left[  X_{1},X_{6}\right]  =\left(

\lambda+2\right)  X_{1}$ & \\

& $\left[  X_{2},X_{6}\right]  =\left(  \lambda+1\right)  X_{2}$ & $\left[

X_{3},X_{6}\right]  =\lambda X_{3}$ & \\

& $\left[  X_{4},X_{6}\right]  =2X_{4}$ & $\left[  X_{5},X_{6}\right]  =X_{5}$%

& \\

$\frak{g}_{6,94}  $ & $\left[  X_{3},X_{4}\right]

=X_{1}$ & $\left[  X_{2},X_{5}\right]  =X_{1}$ & $I_{1}=x_{1}$\\

$\lambda+2=0$& $\left[  X_{3},X_{5}\right]  =X_{2}$ & $\left[  X_{1},X_{6}\right]  =0$ &

$I_{2}=\frac{x_{1}^{2}x_{6}+x_{1}x_{2}x_{5}+2x_{1}x_{3}x_{4}-x_{2}^{2}x_{4}}{2}$\\

& $\left[  X_{2},X_{6}\right]  =-X_{2}$ & $\left[  X_{3},X_{6}\right]

=-2X_{3}$ & \\
\hline

\end{tabular}

\end{table}

\begin{table}[h]

\caption{Nilradical isomorphic to $\mathfrak{g}_{5,6}$ (cont.)}

\begin{tabular}
[l]{llll}%

Name & Brackets & & Invariants\\\hline

$\frak{g}_{6,99}$ & $\left[  X_{3},X_{4}\right]  =X_{1}$ & $\left[

X_{2},X_{5}\right]  =X_{1}$ & none\\

& $\left[  X_{3},X_{5}\right]  =X_{2}$ & $\left[  X_{4},X_{5}\right]  =X_{3}$%

& \\

& $\left[  X_{1},X_{6}\right]  =5X_{1}$ & $\left[  X_{2},X_{6}\right]

=4X_{2}$ & \\

& $\left[  X_{3},X_{6}\right]  =3X_{3}$ & $\left[  X_{4},X_{6}\right]

=2X_{4}$ & \\

& $\left[  X_{5},X_{6}\right]  =X_{5}$ &  & \\\hline

\end{tabular}

\end{table}

\begin{acknowledgement}

The authors wishes to express his gratitute to P. Turkowski for providing him with copies of the articles \cite{Tu1} and \cite{Tu2}, as well as to the Ramon Areces Foundation for support during this research.

\end{acknowledgement}

\end{document}